\documentclass[10pt]{article}
\usepackage{amsmath, amscd, amsthm, amssymb, amsfonts}
\usepackage{verbatim, array, longtable}
\title{Hilbert schemes, Hecke algebras and the Calogero-Sutherland system}
\author{K. Costello and I. Grojnowski}

\newcommand{\counit}{\epsilon}

\newcommand{\disc}{\delta_{\mathrm{disc}}}
\newcommand\rhotw{\tilde{\rho}}
\newcommand{\Del}{\nabla}

\newcommand{\hg}[1]{H_\Gamma^{\otimes {#1}}}
\newcommand{\pro}{\mathrm{pro}}
\newcommand{\alg}{\mathrm{CAlg}}

\newcommand{\iso}{\cong}
\newcommand{\C}{\mathbb C}
\newcommand{\cald}{\mathcal D}
\newcommand{\N}{\mathbb N}
\newcommand{\rarr}{\rightarrow}
\newcommand{\larr}{\leftarrow}
\newcommand{\norm}[1]{\left| #1 \right|}
\newcommand{\Hilb}[2]{\mathrm{Hilb}^{#1} (#2) }
\newcommand{\Oo}{\mathcal O}
\newcommand{\fock}{\mathcal F}
\newcommand{\Z}{\mathbb Z}
\newcommand{\ip}[2]{\left\langle #1 , #2 \right\rangle} 
\newcommand{\lehn}{\mathfrak L}  
\newcommand{\defeq}{\overset{def}{=}}
\newcommand{\cht}{\mathfrak E} 

\newcommand{\be}{\mathbf e}
\newcommand{\ushecke}{\overline{ \mathcal H} } 
\newcommand{\dhilb}{\mathcal {D}^{hilb}}

\newcommand{\aug}[1]{ #1 _+ } 

\newcommand{\Loc}{\operatorname {Loc}}
\newcommand{\Gr}{\operatorname {Gr}}

\renewcommand{\H}{\fock}

\newcommand{\CS}{\mathcal C\mathcal S}
\newcommand{\oCS}{\overline{\CS}}
\newcommand{\tP}{\tilde{\Gamma}}

\newcommand{\IM}{\mathcal I\mathcal M}

 \DeclareMathOperator{\End}{End}
\DeclareMathOperator{\Prim}{Prim}
\DeclareMathOperator{\Supp}{Supp} 
 \DeclareMathOperator{\Ad}{Ad}
\DeclareMathOperator{\Sym}{Sym}
\DeclareMathOperator{\Hom}{Hom}
\DeclareMathOperator{\Spec}{Spec} 
\DeclareMathOperator{\Der}{Der}
\DeclareMathOperator{\Diff}{Diff}
\DeclareMathOperator{\pdiff}{\mathcal {D}iff}
\DeclareMathOperator{\rank}{rank}

\newtheorem{theorem}{Theorem}[subsection]
\newtheorem{proposition}[theorem]{Proposition}
\newtheorem{definition}[theorem]{Definition}
\newtheorem{lemma}[theorem]{Lemma}

\newtheorem{corollary}[theorem]{Corollary}

\newtheorem{remark}[theorem]{Remark}
\newtheorem{example}[theorem]{Example}

\numberwithin{equation}{subsection}

\begin{document}
\maketitle

\begin{abstract}
We describe the ring structure of the cohomology of the
Hilbert scheme of points for a smooth  surface $X$. When $X$ is $\C^2$,
this was done in \cite{lehn-sorger-plane,vasserot} by realising this ring 
as a degeneration
of the center of $\C S_n$. When the canonical class $K_X = 0$,
\cite{lehn-sorger} extended this result by defining an algebra 
structure on  $H^*( \{(x,g) \in X^n \times S_n \mid gx = x \})$; 
the $S_n$-invariants of this algebra is the desired ring.

But when $K_X \neq 0$ it seems
no such algebra can exist.
A completely different approach is needed.

Instead we recast this problem as a question of finding integrals 
of motion for a Hamiltonian which describes ``intersection with the boundary''.
To do so, we use the identification of the cohomology
of the Hilbert schemes with a Fock space modelled on the lattice $H(X)$
\cite{groj,nak}. With this identification, Lehn computed the operator 
of intersection with the boundary. This is essentially the 
Calogero-Sutherland Hamiltonian. 

We then solve the problem of finding integrals of motion by using the 
Dunkl-Cherednik operators to find an explicit commuting family.


We provide two characterizations of 
the Hilbert scheme multiplication operators; the first as
an algebra of operators that can be 
inductively built from functions
and the CS Hamiltonian, and the second in terms of
the centralizer of the CS Hamiltonian inside an appropriate ring 
of differential 
operators.

\end{abstract}

\section{Introduction}

This paper describes the ring structure on the cohomology of the 
Hilbert scheme of points for a smooth 
algebraic surface. The main new
idea in our solution of this problem is to consider the problem as
a question in integrable systems. 
Having done this, we must also {\it solve} 
the integrable
system, i.e{} produce 
integrals of motion.
We do this 
by defining a variant of the 
Cherednik Hecke algebra depending on a finite
dimensional algebra $H$. We then give several characterizations
of the family of commuting differential operators this constructs.

Let us describe the problem we solve, and previous work on it.

\medskip
Let $X$ be an smooth algebraic surface, and $\Hilb{n} X$ the
$2n$-dimensional smooth algebraic variety parameterising length $n$
subschemes of $X$---the $n$'th Hilbert scheme of points for $X$.

The papers \cite{groj, nak}
discovered that the correspondences between Hilbert schemes given by adding points constrained
to lie on cycles in $X$  organize into an action of an infinite dimensional Lie algebra. 
This gives
a canonical identification
$$ \fock =\Sym xH[x] \iso \oplus_n H^*(\Hilb{n} X) $$
as modules for a Heisenberg Lie algebra, and hence in particular an isomorphism of graded 
vector spaces.
(We always work in the category of super vector spaces, so that 
$\Sym V$ means the symmetric algebra of $V^{ev}$ tensored with
the exterior algebra of $V^{odd}$). We write $\fock = \oplus_n \fock^n$, where $\fock^n$ is the $n$'th-eigenspace of the energy operator $\partial = x\frac{\partial}{\partial x}$.

\medskip
The ring structure on cohomology now induces an additional map 
$\fock \otimes \fock \to \fock$,
and the problem that this paper addresses is to describe this map---to determine
the cohomology ring of the Hilbert scheme, in the Fock space coordinates.

\medskip
This problem was first considered by Lehn \cite{lehn}, who also made significant
progress to its solution. Each Hilbert scheme $\Hilb{n} X$ contains a natural divisor,
the complement of the loci of distinct points, and hence a class in $H^2( \Hilb{n} X)$.
Cup product
with this class, summed over all $n$, is thus an operator $\lehn : \fock \to \fock$,
and  $\cite{lehn}$ identifies this operator explicitly. In the Fock space coordinates it
can be written as a 2nd order differential operator depending on the Frobenius algebra 
structure of $H(X)$ and on the canonical class $K \in H^2(X)$.

\medskip
The operator $\lehn$ is naturally associated 
with the first Chern class of the cotangent bundle of ${X}$. 
In a similar manner any cohomology class on $X$ inductively determines
a set of operators on $\fock$, and using the results of \cite{lehn},
Li, Quin and Wang show in a series of papers \cite{lqw, lqw2} that
the subring these operators generate is the ring of multiplication operators on
the Hilbert scheme. After \cite{groj, lehn}, no new algebraic geometry
is needed---these results determine the ring implicitly.

However, 
though much progress
was made in in the above papers,
the relations that these operators satisfy is not clear. In particular, there is no
construction of the operators without already knowing the Hilbert scheme exists,
and no presentation of the cohomology ring.


When $K = 0$, the situation was satisfactorily resolved in \cite{lehn-sorger}:
they construct an algebra
$A_n$ whose underlying vector space is  $H( \{(x,g) \in X^n \times S_n \mid gx = x \})$
and show that its $S_n$-invariants are $H(\Hilb{n}X)$.
This algebra is  the
Chen-Ruan orbifold cohomology ring of $X^n/S_n$, and provides evidence for
a conjecture of Ruan on ring structures on crepant resolutions.

Unfortunately,  this approach fails when $K \neq 0$. 
It seems that
there is no 
flat family of algebras $A_n(K)$ equipped with an $S_n$ action such that $A_n(0) = A_n$ and
whose generic point has $A_n(K)^{S_n} = H(\Hilb{n}{X})$.

New ideas are needed. 

\bigskip
Instead, we 
begin with the following key 
observation.  Cup product with an element $a \in \prod_n \fock^n$ defines 
a linear
map on $\fock$, and, for reasonable $a$, a differential operator on $\fock$. 
This
differential operator commutes with the operator
$\lehn$---after all, cohomology is a (super)-commutative ring.
Hence there are as least as many differential operators which commute with $\lehn$ 
as the dimension of the space $\fock$.

Thinking of $\lehn$ as a Hamiltonian, {\it this is precisely saying $\lehn$
defines an integrable system. It
rephrases the problem of determining the ring structure on $\fock$ 
as the problem of finding integrals of motion}---determining the centralizer
of $\lehn$ in an appropriate ring of differential operators on $\fock$, 
and showing that 
the centralizer is commutative and in bijection with $\fock$.

However, as noticed by 
\cite{bouwk,frenkel_wang,matsuo}, $\lehn$ is a variant of a well known operator---it is a
version of the bosonised Calogero-Sutherland
Hamiltonian, familiar from exactly solvable models of statistical mechanics \cite{amos}. 
So the problem of describing the the ring structure on the Hilbert scheme becomes
an algebraic problem, that of solving a deformed version of the Calogero-Sutherland system.

\medskip
We do that in this paper. We begin
by defining a ring of {\it continuous} differential operators on the Fock space,
$\pdiff_{H,K'} \fock$.
This depends on both the ring structure on $H(X)$ 
and on an additional parameter 
$K' \in H^{ev}(X)$. Section 2 of this paper is devoted to its properties. It is a limit
of finitely generated algebras, and    much smaller than the ring of all differential
operators. 

We also define an operator $\lehn_{K'} \in \pdiff_{H,K'} \fock$;
when $K = K'$ this is the operator $\lehn$.
For generic $K'$, define 
$\IM_{K'}$ to be the centralizer in $\pdiff_{H,K'} \fock$ of $\lehn_{K'}$ 
and of the energy operator $\partial$.
As $\IM_{K'}$ commutes with  $\partial$, for each $n$
it projects to give a subalgebra of $\End \fock^n$. Degenerating $K'$ to $K$,
we get an algebra $\IM_K$ for any $K$.

Then our first description of the ring structure is 
\begin{theorem} $\IM_K$ is a {\em commutative} algebra, 
independent of any choices,
and for any $n$ the image of $\IM_K$ in  $\End \fock^n$ is precisely the algebra
of multiplication operators $H(\Hilb{n}X) \subseteq \End \fock^n$.
\end{theorem}

In order to prove this theorem, we give an explicit construction of $\IM_{K'}$
as a  polynomial algebra. 
We begin by rewriting
the Calogero-Sutherland Hamiltonian $\lehn$ as an inverse  limit 
of $r$-particle generalized
Calogero-Sutherland Hamiltonians.  (To do this precisely forces us
to work with augmented or non-unital algebras).

We then show each $r$-particle system is 
integrable, following the methods of Dunkl, Cherednik and
Opdam. This consists of writing an explicit polynomial algebra of commuting 
difference-differential operators.
The symmetric group acts on these operators, and the invariants act purely as differential 
operators. The quadratic invariant is the CS Hamiltonian, and these operators are the entire 
centralizer of the Hamiltonian. Put differently, we define a version of the degenerate Cherednik
double affine Hecke algebra depending on a Frobenius algbera $H= H(X)$, and its spherical
subalgebra is the desired algebra of integrals of motion.

This gives an explicit description of the integrals of motion for each finite $r$, as well
as showing there are exactly enough of them. Taking the limit over $r$, we get a description
of $\IM_K$ in terms of degenerating families of Dunkl-Cherednik operators.
In particular
we recover the previous theorem.

Section 4 is a description of the generalized Cherednik algebra, and  its
behaviour in the inverse system and flatness properties with respect to twisting.
This enables us to solve the Calogero-Sutherland Hamiltonian $\lehn$ attached to any Frobenius
algebra $H$ and parameter $K$, without needing the existence of the Hilbert scheme to produce
``enough'' commuting operators.

In order to compare the integrals of motion we construct with the
algebra of operators constructed in  \cite{lehn, lqw, lqw2}, we need to 
characterize this latter algebra. It turns out that this has
a simple description as differential operators built out of $\lehn$ 
and the (usual) algebra structure on the Fock space $\fock$. We encode
this in the notion of {\it locality} in section \ref{loc_sec},
and in section 7 we prove that the algebra of differential operators local 
for $\lehn$ is just the desired algebra of multiplication maps.
In section 5 we show that $\IM_K$ are the local operators for $\lehn$
for any $H$ and $K$, proving the main theorem.

\medskip

It seems this is the first occurrence of integrable systems in this manner in algebraic geometry.

\medskip

\paragraph{Acknowledgements.}
We would like to thank Peter Bouwknegt and Ruth Corran for helpful
conversations, and Manfred Lehn for his interest in this work. We
would also like to thank MSRI, where this work was started. K.C.{}
would like to thank the Cecil King Memorial Foundation for
generous financial support.

An earlier 
version of this document has been circulating 
since January 2003, and both authors have given many talks on the
contents since December 2002.

\subsection{Notation.}

All vector spaces are super vector spaces, that is they are
$\Z/2$-graded. By a structure of a commutative algebra on a super
vector space, we mean a super-commutative algebra. Lie brackets 
are taken in the super sense,
that is
\begin{equation*}
[\alpha, \beta] = \alpha \beta - (-1)^{\norm \alpha \norm\beta} \beta \alpha
\end{equation*}
for homogeneous $\alpha$ and $\beta$. Endomorphism algebras of
super vector spaces have a natural $\Z/2$-grading.
If $A$ is a $\Z/2$-graded algebra, we write $A^\times$ for the subgroup of
{\it even} invertible elements.


\medskip
If $X_1, \dots, X_n$ are subsets of an algebra $A$ we write
$\langle X_1, \dots, X_n\rangle$ for the subalgebra of $A$ that they generate.
\medskip

We will often use filtered objects. A filtration on an 
object $X$ in an abelian category is a sequence of subobjects
$F^iX$, $i \in \Z$, such that $F^i X \subseteq F^j X \subseteq X$ if $i \leq j$.
It is {\it exhaustive} if $\cap_i F^i X = 0$ and $\cup_i F^i X = X$,
and non-negatively graded if $F^iX = 0 $ for $i <0 $. Define
$F^\infty X = \cup_i F^iX$ and $F^{-\infty}X = \cap_i F^i X$.
Then there is a functor from filtered objects to exhaustive filtered objects
which sends $X$ to $F^\infty X/F^{-\infty }X$. 

We write $\Gr^i X = F^i X / F^{i-1}X$, and if $Y \subseteq X$ we
give $Y$ the induced filtration $F^i Y = Y \cap F^i X$, so 
$\Gr Y \subseteq \Gr X$.


\section{Differential operators on Fock space }
Let $A$ be a commutative ring. In this section we use the ring
structure on $A$ to define a sequence of ideals in the Fock space
$\Sym A$; i.e.{} to define a topology on $\Sym A$. We then define
the ring of continuous differential operators on $\Sym A$ 
in the standard way. 

In fact, we use the equivalent language of pro-rings. The main result
of this section is
proposition \ref{what_diff}, which identifies $\pdiff \Sym A$.

\subsection{Non-unital algebras} 
We will need to work with augmented algebras. Equivalently, with non-unital
algebras. 
\medskip

Let $\alg_U$ denote the category of unital commutative associative algebras,
$\alg_N$ denote the category of (possibly) non-unital commutative 
associative algebras. The forgetful functor $\alg_U \to \alg_N$
has a left adjoint ${}_+ : \alg_N \to \alg_U$,
$$ \Hom_{\alg_U}(A_+,B) = \Hom_{\alg_N}(A,B) $$
for $A \in \alg_N$, $B \in \alg_U$,
where the underlying vector space of $A_+$ is $\C \times A$,
and multiplication is defined by
$$ (\lambda, u) \cdot (\mu, v) = ( \lambda \mu, \lambda v + \mu u + uv).$$
We write $(1,0) = 1$, and note there are two exact sequences of non-unital
aglebras
$$ 0 \to A \to A_+ \to \C \to 0, \text{ and} $$
$$ 0 \to \C \to A_+ \to A \to 0. $$
If $A \in \alg_U$, these exact sequences split, and there is an isomorphism
of unital algebras $A_+ \iso \C \times A$, 
$(\lambda,u) \mapsto (\lambda, \lambda\cdot 1_A + u)$.
The essential image of $+$ is the category of augmented unital algebras.


Both $\alg_U$ and $\alg_N$ admit coproducts, denoted $\otimes$,
$\otimes_+$ respectively. We have, if $A, B \in \alg_N$ that
$$ (A \otimes_+ B)_+ \iso A_+ \otimes B_+. $$
Note that the underlying vector space of $A\otimes_+ B$ is {\it not}
the tensor product of  the underlying vector spaces of $A$ and $B$.
There is also an action of 
$\alg_U$ on $\alg_N$, denoted $\otimes$; we have
$ (A \otimes B)_+ = A \otimes B_+$, if $A \in \alg_U$, $B \in \alg_N$.

(Recall that the coproduct of $A$ and $B$ is the object in $\alg$
equipped with a morphism from $A$ and a morphism from $B$;
the images of these morphisms generate $A\otimes B$ and the only relation
imposed is the images commute. In $\alg$, the maps
$A \to A \otimes B$, $B \to  A \otimes B$ are injective, 
as the above description shows.)

We will write $(\alg, \otimes)$ to mean either one of 
the categories $(\alg_U, \otimes)$, $(\alg_N, \otimes_+)$ when
no confusion is possible. 

\subsection{Notation}
Let $A \in \alg$ ( = $(\alg_U,\otimes)$ or $(\alg_N, \otimes_+)$).
Write $a \mapsto a_i$ for the $i$'th map $A \to A^{\otimes n}$,
$1 \leq i \leq n$ (recall the definition of $A^{\otimes n}$ as a coproduct).
If $A \in \alg_U$, then 
$a_i = 1\otimes \dots 1\otimes a \otimes 1 \dots \otimes 1$.
More generally, if  $I \subseteq \{1,\dots,n\}$ we have embeddings
$A^{\otimes \# I} \to A^{\otimes n}$, $ \gamma \mapsto \gamma_I$.
For example, given an element $\Delta \in A\otimes A $ we get elements
$\Delta_{ij} \in A^{\otimes n}$.

Similarly, if $\partial \in \Der(A,A)$ is a derivation of $A$, we write
$$ \partial_i : A^{\otimes n} \to A^{\otimes n} $$
for the derivation of $A^{\otimes n}$ which satisfies $\partial_i(a_j) = 0$
if $i \neq j$ and $\partial_i(a_i) = (\partial a)_i$.

\subsection{Infinite symmetric products and symmetric algebras}
\label{inf_sym_prod}

For $A \in \alg_N$, let us denote the $S_n$-invariants of 
$A^{\otimes_+^n}$ by $\Sym^n_+A$. This 
is again in $\alg_N$. We have 
$(\Sym^n_+A)_+  = \Sym^n(A_+)$, where $\Sym^n $ refers to the usual
symmetric power (in either vector spaces or unital algebras). 
The filtration $0 \to \C \to A_+ \to A\to 0$ induces a filtration 
$F^\cdot$
of $\Sym^n A_+$ with  $\Gr_F^a\Sym^n A_+ = \Sym^a A$ if $a \leq n$; i.e.{}
there is a natural inclusion of vector spaces $\Sym^a A_+ \hookrightarrow \Sym^n A_+$
with image $F^a \Sym^n A_+$. This is not an algebra map; instead
each
$\Sym^n A_+$  is a filtered algebra:
$F^a \Sym^n A_+ \cdot F^b \Sym^n A_+ \subseteq F^{a+b} \Sym^{n} A_+$.
We have an isomorphism of vector spaces
$$ F^a \Sym^n A_+ \iso \bigoplus_{0 \le r \le \min(a,n)} \Sym^r A.$$

Dually, the augmentation map $0 \to  A \to A_+ \to \C \to 0$ induces
surjective algebra maps 
$$ p_{nm} : \Sym^{n} A_+ \to \Sym^m A_+, \text{ if $n \geq m$}, $$
and  $p_{nm}(F^a\Sym^{n} A_+) = F^a \Sym^m A_+$ if $a \leq m$.
So  $\Sym^n A_+$ form an inverse system of filtered algebras
in which each filtered piece eventually stabilises.

\begin{definition} \label{def_infinite_sym}
Let
\begin{equation*}
\Sym^{\infty} \aug{A} \defeq F^\infty\lim_{\larr} \Sym^n \aug{A}
\end{equation*}
be the inverse limit taken in the category of filtered algebras.
\end{definition}

We have
$
F^k \Sym^{\infty} \aug{A} \iso F^k \Sym^n \aug {A}
    = \bigoplus_{r \le k} \Sym^r A
$
for any $n \ge k$. 
\medskip
\begin{example}
Fix $t \in \C$, and take $A = \C$ with multiplication $a*b = tab$.
Then on identifying  $\Sym^{\infty} \aug{A} = \C[x]$ as vector spaces,
above, the coefficient of $x^a$ in $x * \cdots *x $ is $t^{n-a}$ times
the number of ways of partitioning $n$ into $a$ non-empty subsets.
For example
$x*x*x*x = x^4 + 6tx^3 + 7t^2x^2 + t^3x$.
\end{example}

Recall that for any vector space $V$,
the vector space $\Sym V = \oplus \Sym^n V $ is the free commutative 
unital algebra
generated by $V$. In particular $\Sym A$ is an algebra 
(non-negatively graded and hence exhaustively filtered).

The natural map $A \to A_+ =  F^1\Sym^nA_+ \to \Sym^nA_+$ lifts to a map 
$A \to  \Sym^\infty A_+$, and hence to an algebra homomorphism
$\Sym A  \to  \Sym^\infty \aug{A}$, which we denote $P$.
Write $P_n : \Sym A \to \Sym^nA_+$ for the quotient maps. 
The sequence of ideals $\ker P_n$ endow $\Sym A$ with 
a seperated topology.

\begin{lemma}\label{surj}
The  morphism
$$ P : \Sym A  \to  \Sym^\infty \aug{A}$$
is an isomorphism of filtered algebras. 
\label{symiso}
\end{lemma}
\begin{proof}
To show surjectivity, it suffices to show that $\Sym^n A_+$
is generated by $A_+ = F^1\Sym^n A_+$. When $r \leq n$ write 
$$\delta_r : A^{\otimes r} \to \Sym^r A \to F^r\Sym^nA_+$$
for the composite of the symmetrisation map 
$a_1\otimes\dots\otimes a_r \mapsto \sum_{w \in S_r} 
a_{w1}\otimes\dots\otimes a_{wr}$ with the inclusion map
(which is itself a symmetrisation). For $r > n$ put $\delta_r = 0$.
Then for $a \in A$,  $\delta_1(a) = P(a)$, and the elements 
$1$, $\delta_r(a_1,\dots,a_r)$ for  $r >0$, $a_i \in A$ span $\Sym^nA_+$.
We have
\begin{equation}\label{symrel}
\delta_r(a_1,\dots,a_r) P(x) = \delta_{r+1}(a_1,\dots,a_r,x) + 
\sum_{1 \leq i \leq r}\delta_r(a_1,\dots,a_ix,\dots,a_r).
\end{equation}
For example, $P(a)P(b) -P(ab) = \delta_2(a,b)$.
From this it is clear that $A$ generates $\Sym^n A_+$. 
Injectivity is obvious.
\end{proof}

\begin{remark}
If $A \in \alg_U$, then the isomorphism $A_+ \iso \C \times A$ induces
an isomorphism of unital algebras
$\Sym^n A_+ \iso \C \times A \times \dots \times \Sym^n A$.
This isomorphism is compatible with the projection maps $p_{nm}$, but not the
filtration.
\end{remark}

\subsection{Differential operators}
In this section we define and identify differential operators on 
the filtered pro-ring $\Sym A$.

\medskip
Recall that a pro-object in a category $\cal C$ is a
functor from $\underline{\N}$ to $\cal C$,
where  $\underline{\N}$ is the category with objects $\{0,1,2,\dots\}$
and a single morphism from $n$ to $m$ if $n \geq m$.
We write  $(X_n, p_{nm} : X_n \to X_m)$ for such a pro-object.
If $(X_n)$, $(Y_n)$ are pro-objects, then
$$ \Hom_{\pro}( (X_n), (Y_n)) \defeq 
\lim_{m \larr } \lim_{\rarr n} \Hom(X_n, Y_m).$$


Furthermore, if $X_n$ and $Y_n$ are filtered pro-objects
then we define
$$ F^i  \Hom(X_n, Y_m) = \{ \varphi \in  \Hom(X_n, Y_m) \mid
\varphi(F^aX_n) \subseteq F^{a+i}Y_m \text{ for all } a \} $$
and put 
$ F^i
 \Hom_\pro( (X_n), (Y_n)) \defeq 
\lim_{m \larr } \lim_{\rarr n} F^i \Hom(X_n, Y_m)$. Now define the
``filtered'' $\Hom$ as  
$$  F^\infty\Hom_\pro( (X_n), (Y_n)) \defeq 
\bigcup_{i \in \Z} F^i  \Hom_\pro( (X_n), (Y_n) ). $$
When each $X_n$ and $Y_n$ are exhaustively filtered this is the correct object;
in general quotient this by $F^{-\infty}$.

\medskip
Recall that if $R$ is a ring, and $M$ an $R$-bimodule, then the
differential part of $M$, denoted $\Diff^{\cdot}M = \bigcup \Diff^s M$
is defined inductively by
$\Diff^s(M)= 0$ if $s < 0$, and 
$\Diff^s(M)= \{m \in M \mid rm -mr \in 
\Diff^{s-1}M \text{ for all } r \in R\}$ if $s \geq 0$.
We write $\Diff^s_R(M)$ if the ring $R$ is not clear from the context.

For example, if $R \in \alg_U$, then
$\Diff \End_\C(R)$ is the ring of differential operators on 
$\Spec R$. When $R$ is not noetherian this is too big a ring.

\medskip
Now let $(A_n, p_{nm})$ be a  commutative pro-ring. Then for $n \geq m$
$\Hom(A_n, A_m)$ is an $A_n$-bimodule, and we may take its differential part.
Suppose that the maps $p_{nm} : A_n \to A_m $ are surjective.
Then for $n \geq m$, the maps
$$ \Hom(A_n,A_n) \to  \Hom(A_n,A_m) \hookleftarrow  \Hom(A_m,A_m) $$
induce maps
$$ \Diff^s\Hom(A_n,A_n) \to  \Diff^s\Hom(A_n,A_m) \hookleftarrow 
 \Diff^s\Hom(A_m,A_m). $$
We have $\Diff^0\Hom(A_m,A_m) \iso  \Diff^0\Hom(A_n,A_m) \iso A_m$,
and so
$$  \Diff^0\End(\lim_{\larr} A_n) = \lim_{\larr} A_n. $$
More generally,
$$\lim_{\larr}\lim_{\rarr}\Diff^s\Hom(A_n, A_m) \hookrightarrow
\Diff^s\End(\lim_{\larr} A_n), $$
but we need not have equality.

Define
$$\Diff^s_\pro \End((A_n)) = \lim_{\larr}\lim_{\rarr}\Diff^s\Hom(A_n, A_m),$$
and $\Diff^._\pro \End((A_n)) = \bigcup_s\Diff^s_\pro \End((A_n))$. 
\medskip

If $(A_n, p_{nm})$ is filtered, then the filtration on $A_n$
induces a filtration on $\Hom(A_n,A_m)$, and thus on 
$\Diff^s \Hom(A_n, A_m)$ and $\Diff^s_\pro\End(( A_n))$.
As always, we redefine the ``filtered'' differential operators to
be the exhaustive part $F^\infty\Diff^s_\pro\End(( A_n)) $.

\medskip
Now, if $(d_n : A_n \to A_n, n \geq 0)$ is such that 
$d_m p_{nm} = p_{nm}d_n$ when $n \geq m$,
then $d_n$ defines an element in $\End(\lim_{\larr} A_n)$.
The set of all such elements form a subalgebra.
Hence if we define
$$ \widetilde{\Diff}^s_\pro\End((A_n))
\defeq \{ (d_n \in \Diff^s \Hom(A_n, A_n),  n \geq 0) \mid
d_m p_{nm} = p_{nm}d_n \} \subseteq \Diff^s_\pro \End((A_n))$$
then
$ \widetilde{\Diff}^._\pro\End((A_n)) $ is a subalgebra of 
$\Diff^._\pro \End((A_n))$. 
We have equality for $s=0$, but not in general.

\begin{example}
$\widetilde{\Diff}^1_\pro  \End ((\C[x]/x^n)) =
x\C[[x]] \subset \C[[x]] =  \Diff^1_\pro \End((\C[x]/x^n)) =
  \Diff^1 \End ( \lim_{\larr}\C[x]/x^n). $
\end{example}

Let $A \in \alg_N$, and consider the pro-object $(\Sym^nA_+, p_{nm})$.
\begin{proposition}
Let $A \in \alg_U$, and $n \geq m \geq 1$. 
If $d \in \Diff^s\Hom(\Sym^nA_+, \Sym^mA_+)$, then
there exists a unique $\tilde{d} \in \Diff^s\Hom(\Sym^mA_+, \Sym^mA_+)$
such that $d = \tilde{d}\circ p_{nm}$.
\end{proposition}
\begin{proof}
We induct on $s$. The case $s = 0$ is clear.
Let $d \in \Diff^s\Hom(\Sym^nA_+, \Sym^mA_+)$. We show that
$d\delta_r(a_1,\dots,a_r) = 0$ for all $r > m$ by descending induction
on $r$. (Our notation is that of lemma \ref{symiso}.) For $r > n$ this is 
clear. Apply $d$ to equation \ref{symrel}, to get
$$
d\delta_r(a_1,\dots,a_r). P(x) + \delta_r(a_1,\dots,a_r). dP(x)
+ d'_x\delta_r(a_1,\dots,a_r) 
 = \sum_{1 \leq i \leq r}
d\delta_r(a_1,\dots,a_ix,\dots,a_r) + d\delta_{r+1}(a_1,\dots,a_r,x),
$$
where $d'_x \in \Diff^{s-1}\Hom(\Sym^nA_+, \Sym^mA_+)$. By induction on $s$,
$d_x' = \tilde{d}_x'\circ p_{nm}$, and $p_{nm}\delta_r(a_1,\dots,a_r) = 0$
for $r > n$. So $d'_x\delta_r = 0$. Also $d\delta_{r+1} = 0 $, by induction
on $r$. Now put $x = 1_A$, to get
$$ (P(1_A) -r ). d\delta_r(a_1,\dots,a_r) = 0.$$
Hence to finish we need only show that $P(1_A) - r $ is not a zero divisor
on $\Sym^mA_+$ for $r > m$. But $P(1_A)$ acts on 
$\oplus_{i \geq k} \Sym^iA / \oplus_{i > k} \Sym^iA 
= \ker(p_{m,k-1})/ \ker(p_{m,k-2})$
as multiplication by $k$, so this is clear.
\end{proof}

The above proposition shows that 
$\Diff^1_\pro\Hom((\Sym^nA_+), A_+) \iso \Diff^1\Hom(A_+, A_+)$.
In contrast with this,
$\Diff^1(\lim_{\larr}\Sym^nA_+, A_+) \iso \Hom(A, A_+)$
(a much bigger space!).

\begin{corollary} \label{difffock}
 $ \widetilde{\Diff}^._\pro\End ((\Sym^nA_+)) 
= \Diff^._\pro\End ((\Sym^nA_+))$. 
\end{corollary}

If $A \in \alg_N$ 
we write $\pdiff \Sym A \defeq F^\infty \widetilde{\Diff}^._\pro
\End ((\Sym^nA_+))$
from now on, and refer to these as {\it the} differential operators
on the Fock space $\Sym A$. As we have seen, they depend on the 
ring structure of $A$. 
The following sequence of propositions describe 
$\pdiff \Sym A$. 

\begin{corollary} Let $A \in \alg_N$. 
i) $F^l\pdiff^s \Sym A = 0$ if $l < 0$.

ii) Let $d \in F^l\Diff^s\End (\Sym^\infty A_+)$.
Then $d \in F^l\pdiff^s \Sym A$ if and only if 
$d \delta_r(a_1,\dots,a_r) \in \ker P_{r-1} $ for all $r$ and all $a_i \in A$.
\end{corollary}

\begin{proposition} \label{nondiff} Let $A \in \alg_U$.
For all $s \geq 0$ there  is a linear map 
$ P^s  : \Diff^s\End_\C(A) \to \pdiff^s \Sym A$ which is uniquely
characterized by requiring that $P^s(D). 1 = P(D(1_A))$, 
and for $s \geq 0$ and all $a \in A$
$$ [P^s(D),P(a)] = P^{s-1}([D,a]). $$ 
In particular, $P^0 = P : A \to \pdiff^0 \Sym A = \Sym A$.
Furthermore, if $D(1_A) = 0$, then $P^s(D) \in F^0\pdiff^s \Sym A$.
\end{proposition}
\begin{proof}
Begin by defining, for $1 \leq i \leq n$, a map
$\Diff^s\End_\C(A) \to \Diff\End_\C(A^{\otimes_+n})$,
denoted $D \mapsto D_i$, by requiring that $$D_i(1) = (D(1_A))_i,$$
where for $a \in A$, we write $a_i$ for the image of $a$ under the
$i$'th coproduct map $A \to A^{\otimes_+n}$,
and for $\Theta \in A^{\otimes_+n}$
$$ D_i(a_j\Theta) = a_jD_i(\Theta)\quad \text{ if } 
j \neq i, \qquad \text{ and } 
\qquad D_i(a_i\Theta) = a_iD_i(\Theta) + [D,a]_i (\Theta). $$
It is immediate that
$(\sum_i D_i). (a_j \Theta) = a_j (\sum_i D_i)(\Theta) + [D,a]_j (\Theta)$,
and that $\sum_i D_i$ is $S_n$-invariant. So $\sum_i D_i$ restricts to
give an element, call it  $P^s_n(D)$, of $\Diff^s\End_\C(\Sym^n A_+)$,
which satisfies the properties of the proposition.

It is clear that $P^0_n = P_n : A \to \Sym^n A_+$,  that for $n \geq m$,
$p_{nm}P^s_n(D) =P^s_m(D) p_{nm} $,  that $P^s(D) \in 
F^1\Diff^s_\pro\End((\Sym^nA_+))$ always; and that if 
$D(1_A) = 0$, then $P^s(D) \in F^0\pdiff^s \Sym A$.
\end{proof}
If $D \in \Diff^s\End_\C(A)$, $D'  \in \Diff^{s'}\End_\C(A)$
then $[D,D'] \in  \Diff^{s+ s' -1 }\End_\C(A)$ and we have
$$ [P(D), P(D') ] = P[D,D']. $$
So $P$ is a map of $\C$-Lie algebras. The enveloping algebra of the source 
consists of ``non-linear differential operators''.

\medskip
Write $\pdiff \Sym^n A_+ = \{ d \in \Diff\End \Sym^nA_+ \mid
d (\ker p_{nm}) \subseteq \ker p_{nm} \text{ for all } m \leq n \}$.
The natural map $p_n : \pdiff \Sym A \to \pdiff \Sym^n A_+$ admits a  section,
defined as follows.

If $D \in \pdiff \Sym^n A_+$, then composing with the symmetrisation map
$\be = \frac{1}{n!}\sum_{w \in S_n}w$, we get a map, also denoted $D$,
$ A^{\otimes_+n} \xrightarrow{\be}  (A^{\otimes_+n})^{S_n}\xrightarrow{D}
 (A^{\otimes_+n})^{S_n} \hookrightarrow  A^{\otimes_+n} $.
Now if $r \geq n$, and $I \subseteq \{1,\dots,r\}$ is a subset such
that $\# I = n$, let 
$D_I :  A^{\otimes_+r} \to A^{\otimes_+r} $
be the map which is $D \otimes Id$, $D$ in the $I$'th places.
Finally, put $\Gamma_n(D)_r = \sum_{I} D_I  $; the sum is over subsets 
$I \subseteq \{1,\dots,r\}$ of size $n$.

It is immediate that $\Gamma_n(D)_r$ is an $S_r$-invariant differential operator;
that $\Gamma_n(D)_n = D$, and that $p_{rk}\Gamma_n(D)_r = \Gamma_n(D)_r$ for 
all $r \geq k \geq n$. Hence we have defined an injection
$$ \Gamma_n :  \pdiff \Sym^n A_+ \to \pdiff \Sym A,$$
such that $p_n \Gamma_n = Id$; in particular $p_n$ is surjective.
Note that $\Gamma_n$ is not an algebra homomorphism, and that the image of
$\Gamma_n$ is in $F^n \pdiff \Sym A$.

\begin{proposition} If $D \in F^r\pdiff^s \Sym A$, then 
$D = \Gamma_{r+s+1} p_{r+s+1}(D)$. In particular,
$$\pdiff \Sym A  = \lim_{\rarr n} \pdiff \Sym^n A_+. $$
\end{proposition}
\begin{proof} A differential operator of order $s$ is determined by
its values on products of $\leq s$ elements. As $\Sym A$ is generated by 
$F^1 \Sym A = A$, an element $d \in \pdiff^s \Sym A$ is determined 
by its values on $F^s\Sym A$. If $d \in  F^r\pdiff^s \Sym A$,
then $d(F^s \Sym A) \subseteq F^{r+s} \Sym A$.
Hence if $d, d' \in F^r\pdiff^s \Sym A$
induce the same map on $\pdiff \Sym^{r+s+1} A_+$ then $d= d'$. Taking
$d' = \Gamma_{rs+1} p_{rs+1}(d)$ we get the proposition.
\end{proof}

To finish the description of $\pdiff \Sym A$, we need only describe
$ \pdiff \Sym^n A_+$. If $A \in \alg_U$, then as 
$\Sym^n A_+ \iso \oplus_{0 \leq i \leq n} \Sym^i A $, 
$\pdiff \Sym^n A_+ \iso \oplus_{0 \leq i \leq n} \Diff \End\Sym^i A $.
Hence $\Diff \End\Sym^n A$ embdeds into $ \pdiff \Sym^n A_+$.
Write $\tP_n :  \Diff \End\Sym^n A \to \pdiff^s \Sym A$ for the 
composite of $\Gamma_n$ with the embedding. For example $\tP_1 = P$ is
the map defined in proposition \ref{nondiff}.

\begin{corollary} \label{what_diff} Let $A \in \alg_U$.
 There is an isomorphism of vector spaces
$$
\oplus_{n \geq 0} \tP_n : \oplus_{n \geq 0} \Diff^s \End \Sym^n A \iso 
\pdiff^s \Sym A.
$$
\end{corollary}

\subsection{Frobenius algebras}
A non-unital {\it Frobenius algebra} $H$ is an algebra $H \in \alg_N$
equipped with a map $\Delta : H \to H \otimes H$ such
that $\Delta$ is an $H$-bimodule map: 
$\Delta(ahb) = a \Delta(h) b$. If $H \in \alg_U$, then $\Delta$ is 
determined by $\Delta(1_H)$ and we call $H$ a {\it weak } Frobenius algebra.
If in addition $\Delta$ has a
counit $\counit : H \to \C$ such that $\counit \otimes Id \circ  \Delta = 
Id \otimes \counit \circ\Delta  = Id$, 
then $H$ 
is a Frobenius algebra in the usual sense.

We assume in the above definitions that all maps are even.

If $H$ is a weak  Frobenius algebra, the element $e = m \Delta(1_H) \in H$ is called the
{\it Euler class} of $H$; here $m : H \otimes H \to H$ denotes the multiplication map.
So if $\Delta(1_H) = \sum_j a_j \otimes b_j$, $ e = \sum_j a_j b_j$.

\subsection{Relative differential operators} 

The following definition is somewhat ad hoc, but will do for our purposes.

Suppose $A = H \otimes R = H_R$, 
where $ H$ is a weak Frobenius algebra, and $R$ is a localisation
of $\C[x]$. (We will only use $\C[x]$ or $\C[x,x^{-1}]$ below).
Let $n \geq 1$, and  $\disc = \prod_{i < j}(x_i - x_j) \in R^{\otimes n}$
be the discriminant.  
Write $\Delta = \Delta(1_H)$.

Now define the {\it relative differential operators}
$$ \Diff_H \End \Sym^n H_R =  
\Big\{ d \in 
 \Diff \End_{H^{\otimes n}}
( H_R^{\otimes n}[\frac{\Delta_{ij}}{x_i - x_j}\mid i < j])
\,\Big|\, d \Sym^n H_R \subseteq \Sym^n H_R \Big\}, $$
where $ \End_{H^{\otimes n}}$ refers to $H^{\otimes n}$-linear 
endomorphisms, and $ H_R^{\otimes n}[\frac{\Delta_{ij}}{x_i - x_j}\mid i < j]$
is the subring of the localisation $ H_R^{\otimes n}[\frac{1}{\disc}]$
generated by $H_R^{\otimes n}$ and the elements $\frac{\Delta_{ij}}{x_i - x_j}$.

 Clearly this is an algebra; if $H = \C$ and $\Delta \neq 0$ it is 
just $\Diff \End \Sym^n R$.
This construction works whether we intepret $H_R$
to be in $(\alg_U, \otimes)$ or  $(\alg_N, \otimes_+)$;
hence repeating the discussion of the previous section we can define
 \begin{align*}  \pdiff_H \Sym H_R  
& = F^\infty \widetilde{\Diff}_{H,\pro} ((\End \Sym^n (H_R)_+)) \\
& \iso \oplus \tP_n  \Diff_H \End \Sym^n H_R.
\end{align*}

\subsection{Twisting} \label{sub_twist}
For $u \in A$, define an algebra endomorphism
$$ \Phi_u : \Sym A \to \Sym A $$
by $P(a) \mapsto P(ua)$ if $a \in A$. 
For simplicity we assume throughout that $u$ is even.
We have $\Phi_u\Phi_v = \Phi_{uv}$, for $u, v \in A$, so if $u \in A^\times$ then
$\Phi_u^{-1} = \Phi_{u^{-1}}$.

\begin{lemma} We have $ \Phi_u( F^r \Sym^\infty A_+)  \subseteq  F^r \Sym^\infty A_+$. More precisely,
$$\Phi_u \delta_r(a_1,\dots,a_r) = \delta_r(ua_1,\dots,ua_r) + \cdots + 
P(h_r(u).a_1\dots a_r) \in F^r \Sym^\infty A_+,$$
where $h_r(u) = u(u-1)\dots(u-r)$.
\end{lemma}
It follows that for  $u\not \in \Z$, $\Phi_u(\ker P_n) \neq \ker P_n$, and moreover
that $\Phi_u$ is not continuous. (It is  continuous if for all $m >1$, there exists
an $n > m$ such that $\Phi_u(\ker P_n) \subseteq \ker P_m$).

We have $\Phi_u \pdiff^0 \Sym A \,\Phi_u^{-1} = \pdiff^0 \Sym A = \Sym A$.
If $\partial \in \Der(A,A)$, then 
$$ [\Phi_u P(\partial) \Phi_u^{-1} , P(a) ] = P(\partial a + \frac{\partial(u^{-1})}{u^{-1}} \cdot a) ,$$
so that if $\partial(u) = 0$, then 
$\Phi_u P(\partial) \Phi_u^{-1} = P(\partial)$.
Hence if $\partial_i \in \Der(A,A) $ and $\partial_i u = 0$ for $i = 1,2$,
$$  [\Phi_u P(\partial_1 \partial_2) \Phi_u^{-1} , P(a) ] =
P(\partial_1 \partial_2(a) + u^{-1} \partial_1(a).\partial_2 + u^{-1} \partial_1. \partial_2(a)), $$
so that  $ \Phi_u P(\partial_1 \partial_2) \Phi_u^{-1} $ is not in $P(\Diff^s \End A)$ or in
$\pdiff^2 \Sym A$.

\begin{lemma} i) Suppose $e \in H$ is nilpotent. Then the map 
$ u \mapsto u - eu^{-1}$ defines a surjection $H^\times \to H^\times$.

ii){\it  (The splitting principle.)}
Suppose $K \in H$ is given. Then there is a Frobenius algebra $\tilde{H}$
containing $H$ as an index $2$ subalgebra, and an element  $u \in \tilde{H}$ such that 
$ u^2 - Ku - e = 0 $ in  $\tilde{H}$.
\end{lemma}
\begin{proof} i) As $e$ is nilpotent, $K = u(1-eu^{-2})$ is invertible.
Conversely, suppose $K \in H^\times$.
Then take $ u = \frac{1}{2}K.(1 + \sqrt{(1 + \frac{4e}{K^2}}).$

ii) Define $\tilde{H} = H[u]/ u^2 - Ku - e$. 
\end{proof}
We do not use part (ii) of the lemma.

Define a {\it degeneration direction} to be a pair of maps $\C \to H$,
$\C^\times \to H$, denoted 
$\lambda \mapsto K_\lambda$, $\lambda \mapsto u_\lambda$ such that for
each $\lambda \neq 0$, $u_\lambda$ is invertible and
$u_\lambda - e u^{-1}_\lambda = K_\lambda$. Write $K = K_0$.
For example, if $e$ is nilpotent and $K \in H$,
then
$(\lambda + K,  \frac{1}{2}K_\lambda(1 + \sqrt{(1 + \frac{4e}{K_\lambda^2}})$
is a degeneration direction; if $e^2 = 0 = eK$, this is just
$(\lambda + K, \lambda + K + \frac{e}{\lambda})$.

\begin{definition} Let $(K_\lambda,u_\lambda)$ be a degeneration direction.
Define
$$ \pdiff_{H,K} \Sym H_R =  \lim_{\lambda \to 0}
\Phi_{u_\lambda}
 \, \pdiff_H \Sym H_R \, \Phi_{u_\lambda}^{-1}. $$
\end{definition}
The notation is abusive; it omits the choice of degeneration direction.
Note that the limit exists as a subalgebra of $\Diff\End\Sym H$ for general 
reasons, but that these general reasons give little control over the
resulting algebra.

\subsection{Locality} \label{loc_sec}

Let $\cald_A \subseteq \Diff \End(A) $ be a subalgebra of differential
operators on $A$ such that $A \subseteq \cald_A$. 
\begin{definition}
Let $B \subseteq \cald_A$.
An algebra $E \subset \cald_A$ is said to be \emph{local} with
respect to $B$, if
\begin{equation*}
  [B,E] = 0 \qquad \text{ and } \qquad
[E,A] \subset \ip A B \cdot E.
\end{equation*}
\end{definition}

Suppose also that $[B,B] = 0$, so that $B \subseteq Z_{\cald_A}(B)$.
If $E$ is such that $[B,E] = 0$, then 
$[E,A] \subset \ip A B  E$ $\Leftrightarrow$
$[E,A]  \subset E \ip A B $ $\Leftrightarrow$
$E \ip A B = \ip A B E $ $\Leftrightarrow$ 
$E \ip A B $ is a subalgebra of $\cald_A$.

\medskip
If $E'$ and $E''$ are local with respect to $B$ 
 then so is $\ip {E'} {E''}$.
Hence there is a maximal 
algebra local with respect to $B$,
which we denote $\Loc_{\cald_A}(B)$. So 
$B \subseteq \Loc_{\cald_A}(B) \subseteq  Z_{\cald_A}(B)$.
If $\ip{A}{\Loc_{\cald_A}(B)} = \cald_A$, then 
$ \Loc_{\cald_A}(B) =  Z_{\cald_A}(B)$.

\medskip

$\Loc(B)$ is filtered. For any
$X \subseteq  \Diff \End(A)$ write $X^s = X \cap  \Diff^s \End(A)$.
Then
$\Loc(B)^{-1}= 0$, and $\Loc(B)^s = \{ d \in \cald^s_A \mid [d,B] = 0 
\text{ and }
[d,A] \subseteq \sum_{0 \leq i \leq s-1} {\ip A B}^i \Loc(B)^{s-1-i} \}$.

\medskip
We will apply these notions to $\pdiff_{H,K} \Sym H_R$ and $\Diff\End \fock$.
\section{The Calogero-Sutherland system}
\label{section_cs}

In this section we define a generalized Calogero-Sutherland type 
Hamiltonian $\lehn$, which will turn out to be integrable.
This is a second order differential operator acting on the Fock space 
$\Sym xH[x]$, where $H$ is a Frobenius algebra. 
In section \ref{section_hs} we will take $H$ to be the cohomology of a smooth
surface, and then the integrals of motion for this Hamiltonian
we will be precisely the operations of multiplication on the 
individual Hilbert schemes. In the next section we describe the 
integrals of motion for arbitrary $H$ in terms of a Cherednik Hecke algebra.

\medskip

Let $H$ be a vector space, $\Gamma =\C[x,x^{-1}]$
and write $H_\Gamma = H \otimes \Gamma = H[x,x^{-1}]$.
Write 
$\partial = \partial_x =  x\frac{\partial}{\partial x}$, an $H$-linear 
derivation of $H_\Gamma$.

Define the {\it Fock space} $\fock(H)$ to be 
$$ \fock(H) = \Sym H[x] / H .  \Sym H[x], $$
so that $\fock(H)$ is a subquotient of $\Sym H_\Gamma$ which
is isomorphic to 
$ \Sym xH[x]$  as a vector space.

Now suppose that $H \in \alg$ is an algebra. Then we define 
first order differential operators 
$P(hx^a \partial) \in  \pdiff^1 \Sym H_\Gamma \subseteq \Diff^1 \End \Sym H_\Gamma $ for
each $hx^a \in H_\Gamma$ as in \ref{nondiff}; i.e.{} by requiring that
$$ [P(hx^a \partial), P(h'x^n) ] = n P(hh' x^{a+n}),
\qquad P(hx^a \partial).1 = 0. $$
It is clear that for $a \geq 0$ these operators descend to give differential
operators on $\fock(H)$. The operator $P(\partial)$ is called
the {\it energy} operator, and its eigenspaces the {\it energy weight spaces}.
If $H$ is finite dimensional then the eigenspaces of $P(\partial)$ on
$\fock(H)$ are finite dimensional; those on $\Sym H_\Gamma$ are not.

Further suppose $H \in \alg_U$ is a weak Frobenius algebra,
so $\Delta = \Delta(1_H) \in H \otimes H$. 
Write  $\Delta = \sum_j a_j \otimes b_j$, $e = \sum_j a_jb_j \in H $ and
define 
$ \Delta_* : H_\Gamma \to \Sym H_\Gamma $
by 
$$ \Delta_* P(hx^n) = \sum_{j; r \in I(n)} P(a_jhx^r)P(b_jx^{n-r}), $$ 
where  
$I(n) = \{ 1,\dots,n\}$ if $n \geq 0$, and
$I(n) = \{n,\dots,-1\}$ if $n \leq -1$.
Notice that $ \Delta_* P(hx^n) \in  |n|.P(ehx^n) + \ker P_1$.

Let $K \in H$. Then we define a 
second order differential operator 
$\lehn = \lehn(H,K) \in \Diff^2 \End \Sym H_\Gamma $
by requiring that $ \lehn. 1 = 0$ and 
$$ [\lehn, P(hx^n) ] = 2n P(hx^n\partial) + n^2P(Khx^n) +
|n| \Delta_* P(hx^n)
$$
Again, it is clear that this descends to give a differential
operator on $\fock(H)$, the {\it Calogero-Sutherland operator}.

We retain this notation for the action of these operators on any 
invariant subquotient of $\Sym H_\Gamma$.

\begin{proposition} \label{ctsCS}  
Write $H_\Gamma^{\otimes 2} = H^{\otimes 2}[x^{\pm 1}, y^{\pm 1}]$,
and let $u \in H^\times$ be any even invertible element.
Then in the notation of section \ref{sub_twist}
$$ \lehn(H,u-eu^{-1}) = 
\Phi_u \, \Big(\tP_1(u\partial^2) + 
\tP_2(\frac{x+y}{x-y} u^{-1} \Delta (\partial_x - \partial_y)) \Big)
\, \Phi_u^{-1}.$$

In particular, $ \lehn(H,K) \in F^1\pdiff^2_{H,K} \Sym H_\Gamma$ 
for all degeneration directions $(K_\lambda,u_\lambda)$.

\end{proposition}
\begin{proof}
Write $\lehn' = \Phi_u^{-1} \lehn(H, u-eu^{-1}) \Phi_u$,
and  $$\lehn'' = 
\tP_1(u\partial^2) + 
\tP_2(\frac{x+y}{x-y} u^{-1} \Delta (\partial_x - \partial_y)).$$

We must show that 
 $\lehn' = \lehn''$.
This is straightforward from the definitions. Begin by observing
\begin{align*}
   [\lehn'' , P(hx^n) ] &=n^2P(uhx^n) +  2nP(uhx^n\partial) + 
n \sum_{i<j}h_i u_i^{-1} \Delta_{ij}
\frac{x_i + x_j}{x_i - x_j} (x_i^n - x_j^n) \\
&= n^2P(uhx^n) + 2nP(uhx^n\partial) + |n| \sum_{i,j; r \in I(n)} 
h_iu_i^{-1} \Delta_{ij}x_i^rx_j^{n-r}  
- n^2 \sum_i h_i u_i^{-1}e_ix_i^n \\
& = 2nP(uhx^n\partial) +  n^2 P((u-eu^{-1})h x^n) + 
|n|\Delta_* P(u^{-1}hx^n).
\intertext{
But we also have
}
 [  \lehn' , P(hx^n) ] & = 2n P(u hx^n \partial) + n^2 P( (u-eu^{-1})hx^n) + 
 |n| \Delta_* P(u^{-1}hx^n) .
\end{align*}

As $\lehn' 1= \lehn'' 1 = 0$, we have $\lehn' = \lehn''$.
Finally, observe that as we have equality on the subset $H^\times \subseteq H$,
it must be that the $\Phi_u \lehn'' \Phi_u^{-1}$
  depends  only on  $u-eu^{-1}$. Hence for
any degeneration direction $(K_\lambda,u_\lambda)$ its limit is just 
$\lehn(H,K) \in F^1\Diff^2$.

\end{proof}

\begin{remark} 
An easy direct computation 
shows that 
$ \lehn'. \big( P(hx^n)P(h'x^m) - P(hh'x^{n+m}) \big) \in \ker P_1$ 
if and only if $u - u^{-1}e = K$. Hence
$ \lehn = \lehn(H,K) \in \Phi_u \pdiff^2 \Sym H_\Gamma \,\Phi_u^{-1} $
if and only if $u - u^{-1}e = K$.
\end{remark}

Define $\IM_K$ to be the centralizer in $\pdiff_{H,K} \Sym H_\Gamma$ of
both $\lehn(H,K)$ and $P(\partial)$,
$$ \IM_K = Z_{\pdiff_{H,K} \Sym H_\Gamma}(\lehn(H,K), P(\partial)) $$
if $K = u-eu^{-1}$ for some invertible $e$, and
in general let $\IM_K = \lim_{\lambda \to 0}\IM_{u_\lambda - e u_\lambda^{-1}}$
for a choice of degeneration direction $(K_\lambda, u_\lambda)$.
As defined $\IM_K$ seems to depend on the choice of degeneration direction. 
The notation is acceptable, because of 
the following theorem, 
which will be proved in section \ref{section_hecke}, and in 
theorem \ref{construct}.

\begin{theorem} 
i) The algebras $\IM_K$ depend only on $K$, as algebras, subalgebras of $\End \Sym H_\Gamma$,
and pro-finitely generated algebras.

ii) $\IM_K$ forms a flat family of algebras as $K$ varies. 

iii) $\IM_K$ is commutative.
\end{theorem}

\section{Hecke algebras}
\label{section_hecke}
In this section we define a variant of the Dunkl-Cherednik
operators. These operators act on $\hg{n}$, where $H$ is a weak 
Frobenius algebra.
As in the usual case these operators commute, and together with the 
group algebra $\C S_n$ form a generalization of the
degenerate affine Hecke algebra of 
type $A_n$. The algebra obtained by including the operations of multiplication
by elements of $\hg{n}$ also closes; this is a generalization of 
the Cherednik degenerate double affine Hecke algebra.

Our exposition follows closely \cite{heck} (which considers the rational case).
See also \cite{chered2,opdam1}. When $H = \C \in \alg_U$, the results 
in sections \ref{sec_notheck}--\ref{sec_prop} are contained in these papers.

In sections \ref{sec_notheck}--\ref{sec_prop} all constructions work
in either of the categories $(\alg_U, \otimes)$ or $(\alg_N, \otimes_+)$.
In section \ref{sec_stable} we must specialise to the non-unital case
in order to get augmentation maps.

\subsection{Notation}
\label{sec_notheck}

Let $A \in \alg$ ( = $(\alg_U,\otimes)$ or $(\alg_N, \otimes_+)$),
and  $\partial \in \Der(A,A)$.
If $l = (l_1,\dots,l_n) \in \Z^n$, write
$$ \partial_l = \sum_{i} l_i\partial_i \in \Der(A^{\otimes n},A^{\otimes n}).$$
We have $\partial_{l+l'} = \partial_l + \partial_{l'}$.

Equip $\Z^n$ with the inner product
$\langle l,l' \rangle = \sum l_il'_i$, and let
$\Phi = \{ l \in \Z^n \mid \langle l,l \rangle = 2 \} = 
\{ \varepsilon_i - \varepsilon_j \mid i \neq j \}$,
where $\varepsilon_i = (0,\dots,0,1,0,\dots,0)$.
Write $\alpha > 0$ to mean $\alpha = \varepsilon_i - \varepsilon_j \in \Phi$ with
$ i < j $. 

If $\alpha = \varepsilon_i - \varepsilon_j$ define 
$$ r_{ij} = r_\alpha : \Z^n \to \Z^n, \qquad
l \mapsto l - \langle \alpha, l\rangle \alpha $$
to be the associated reflection, so that the group generated by the relections
$r_\alpha$ with $\alpha \in \Phi$ is just $S_n$.

Observe that $S_n$ acts on $A^{\otimes n}$, and that for $\partial \in
\Der(A,A)$,
$$ r_\alpha \partial_l r_\alpha^{-1} = \partial_{r_\alpha l}. $$

\subsection{Dunkl-Cherednik operators}
We now fix $A$ to be $H_\Gamma = H \otimes \C[x,x^{-1}]$ where
$ H \in \alg_U$. The results in this section make sense
if we consider $H_\Gamma$ to be in either $(\alg_U,\otimes)$ or
$(\alg_N,\otimes_+)$.
 Our notation is such that $\hg{n}$ is generated by 
the elements $(hx^a)_i = h_i x_i^a$, with $1 \leq i \leq n$, 
$h \in H$, $a \in \Z$. Further suppose that $H$ is a weak Frobenius algebra,
and let $\Delta = \Delta(1_H) \in H \otimes H$.
Finally, fix the $H$-linear derivation 
$\partial = x\frac{\partial}{\partial x} \in \Der(H_\Gamma,H_\Gamma)$.

We write for $ \alpha = \varepsilon_i - \varepsilon_j \in \Phi$,
$e^\alpha = x_ix_j^{-1}$, and $\Delta_\alpha = \Delta_{ij} \in \hg{n}$,
and define
$$ \Del_\alpha = \Del_{ij} : \hg{n} \to \hg{n} \qquad\text{ by }\qquad
\Del_{ij} = \frac{1}{1-x_ix_j^{-1}} \Delta_{ij} (1 - r_{ij}). $$
As $r_{ij} \Delta_{ij} = \Delta_{ji} = \Delta_{ij}$ and
$\Delta_{ij} (h_i - h_j) = 0$, we have
$$ \Del_{ij} \left( (hx^a)_i\, (h'x^b)_j \right) = 
(hh')_i \Delta_{ij} \frac{x_i^a x_j^b - x_i^b x_j^a}{1-x_i x_j^{-1}}, 
\text{ and }$$
$$ \Del_{ij} \left( (hx^a)_k \, \Theta \right )   = (hx^a)_k \, \Del_{ij}(\Theta)
 \qquad\text{ if } \Theta \in \hg{n} \text { and } k \not\in \{i,j\}. $$
For $ u \in H$, write $(u\Del)_{ij} = u_i \Del_{ij} =  u_j \Del_{ij}$.
Finally,
let $\rhotw(u) = \frac{1}{2} \sum_{\alpha > 0}(u\Delta)_\alpha\alpha$.

\begin{definition} Let $l \in \Z^n$, $u \in H^\times$. 
The Dunkl-Cherednik operator with parameter $u$
is defined to be
$$ y_l = y_l(u) = 
\partial_l + \sum_{\alpha > 0} \langle \alpha, l \rangle (u^{-1}\Del)_\alpha
-  \langle \rhotw(u^{-1}), l \rangle 
: \hg{n} \to \hg{n}.$$
\end{definition}
We write $y_i$ instead of $y_{\varepsilon_i}$. Define an action of $S_n$
on the Dunkl operators by ${}^w\!y_l = y_{wl}$.

\begin{proposition}
\label{Dunkcomm}
i) $[y_l, y_{l'}] = 0 $ if $l$, $l' \in \Z^n$.

ii) For $f(y)  \in \C[y_1,\dots, y_n]$, and 
$\alpha = \varepsilon_i - \varepsilon_{i+1}$,
$$ f. r_{\alpha} - r_{\alpha} .{}^{r_{\alpha}}\!f = (u^{-1}\Delta)_{\alpha} \,
\frac{f - {}^{r_{\alpha}}\!f }{y_i - y_{i+1}} $$
as operators from $\hg{n}$ to $\hg{n}$.
\end{proposition}

The most straightforward proof is by direct computation. We recall
some of the main steps.
\begin{lemma} 
$$i) \qquad 
\left[ \partial_l, \Del_\alpha \right] = 
\frac{ \langle \alpha, l \rangle}{1- e^{-\alpha}}
\left( \Delta_\alpha r_\alpha \partial_\alpha - 
\partial_\alpha(e^{-\alpha}). \Del_\alpha \right). $$

$$ ii) \qquad \left[ y_l, y_{l'} \right] = \sum_{\alpha, \beta > 0}
\left( \langle \alpha, l \rangle \langle \beta, l' \rangle -
 \langle \alpha, l' \rangle \langle \beta, l \rangle \right) 
(u^{-1}\Del)_\alpha (u^{-1}\Del)_\beta .$$
\end{lemma}
Explanding out $\Del_\alpha\Del_\beta$ this reduces the problem to 
checking it in the rank two root systems. Arguing as in \cite{heck} 2.2,
the proposition reduces to the identity 
$$ -(1-z) + (1-w) \frac{1-z}{1-z^{-1}} + (1-wz) = 0,$$
which we apply with $z = x_i/x_j$,  $w = x_j/x_k$ and  $i$, $j$, $k$ all 
distinct.
We omit further details.

The proposition can rephrased: Define the {\it {Cherednik algebra}}
(or {\it {degenerate double affine Hecke algebra}}) to 
be the vector space
$$ \ushecke_n = \hg{n} \otimes \C S_n \otimes \C[y_1,\dots, y_n]$$
with the unique algebra structure that makes it a subalgebra of
$\End_{\C}(\hg{n})$ by which $y_l$ act as Dunkl operators, $S_n$ act
as permutations, and $\hg{n}$ acts by multiplication. Namely, we require 
that $\C[y_1,\dots, y_n]$, $\C S_n$ and $\hg{n}$ are subalgebras,
and that the relations of \ref{Dunkcomm}(ii) hold, and that 
$[y_l, \Theta] = y_l(\Theta)$, $\Theta\, w = w \,{}^w\!\Theta$ for 
$\Theta \in \hg{n}$, $l \in \Z^n$, $w \in S_n$.
(For this to make sense, always regard $\C S_n$ and 
$\C[y_1,\dots, y_n]$ as unital algebras when taking tensor product with
$H_\Gamma^{\otimes n}$).

Notice that  $\C S_n \otimes H[y]^{\otimes n}$ forms a subalgebra
of $\ushecke_n$. This algebra has the same relation to the
``degenerate affine Hecke''
algebra of Drinfeld and Lusztig that 
$\ushecke_n$ has to the usual Cherednik algebra.

\subsection{Properties}
\label{sec_prop}

A direct computation shows
\begin{lemma} \label{cascomp}
 $$ i) \qquad  \sum_{1 \leq i \leq n} y_i 
 =  \sum_{1 \leq i \leq n} \partial_i, \text{ and } $$
 $$ ii) \qquad \sum_{1 \leq i \leq n} y_i^2  
 - \langle \rhotw, \rhotw \rangle
= \sum_{1 \leq i \leq n} \partial_i^2 + 
\sum_{\alpha > 0}  \frac{1+e^{-\alpha}}{1-e^{-\alpha}} (u^{-1}\Delta)_\alpha \partial_\alpha
 + \Xi,$$
where $\Xi(\Sym^n H_\Gamma) = 0 $.
\end{lemma}
More generally we have
\begin{lemma} \label{dif_lem_0}
The action of $H[y]^{\otimes n}$
on $\hg{n}$ restricts to an action of $\Sym^n H[y]$ on 
$\Sym^n H_\Gamma$. Furthermore, this action is by  differential
operators which preserve the energy weight spaces, and
 $ \Sym^n H[y] \subseteq \Diff_H \End \Sym^n H_\Gamma$.

\end{lemma}
\begin{proof} By \ref{Dunkcomm}(ii), if $f \in \Sym^n H[y]$ then $wf = fw$ for
all $w \in S_n$, so $f$ preserves $\Sym^n H_\Gamma$. To show this action 
is by differential operators, it suffices to show that $H[y]^{\otimes n}$ is in
$\Diff_{\Sym^n H_\Gamma}\Hom(\Sym^n H_\Gamma,H_\Gamma^{\otimes n})$.
But 
$(hy)_i = (hu\partial)_i \in 
\Diff^1\Hom(\Sym^n H_\Gamma,H_\Gamma^{\otimes n})$
and these
elements generate $H[y]^{\otimes n}$.
This suffices, as it is immediate from the definition of $\Diff_R^sM$
that if $D_1, D_2 \in \End(H_\Gamma^{\otimes n})$ are such that
they restrict to differential operators 
$ D_i \in \Diff^{s_i}\Hom(\Sym^n H_\Gamma,H_\Gamma^{\otimes n})$
then $D_1D_2  \in \Diff^{s_1+s_2}\Hom(\Sym^n H_\Gamma,H_\Gamma^{\otimes n})$.
\end{proof}

We can rephrase the lemma: Write $\C[y_1,\dots, y_n]_{\leq s}$ for
polynomials in $y$ of total degree less than $s$;
the degree of each $y_i$ is 1.
\begin{corollary}  \label{diff_lem} 
The $\Sym^n H_\Gamma$-differential part of 
$\ushecke_n$ is all of $\ushecke_n$. Moreover, this is a filtration
of algebras:
$\Diff^i \ushecke_n. \Diff^j \ushecke_n \subseteq \Diff^{i+j} \ushecke_n$,
with
$\Diff^s_{\Sym^n H_\Gamma}\ushecke_n = 
H_\Gamma^{\otimes n} \otimes\C S_n \otimes \C[y_1,\dots, y_n]_{\leq s}$.

Finally, write $p_i$ for the image of $y_i$ in $\Gr \ushecke_n$. 
Then we have an algebra isomorphism
$$ \Gr \ushecke_n  = (H \otimes \C[x,x^{-1},p])^{\otimes n} \# \C S_n,
$$
where $S_n$ acts on the $n$-fold tensor power in the obvious way.
\end{corollary}

Write $P : \Sym H_\Gamma \otimes \Sym H[y] \to 
\Diff \End \Sym^n H_\Gamma $ for the map
$f(x) \otimes g(y) \mapsto \sum_{1\leq i \leq n} f(x_i)g(y_i)$.
This notation is compatible with all previous maps called $P$.

So lemma \ref{cascomp}(i) shows $P(y) = P(\partial)$,
$[P(hx^ay), P(h'x^m)] = mP(hh'x^{a+m})$
and \ref{cascomp}(ii)
says 
$$ P(y^2) - \langle \rhotw, \rhotw \rangle
 = \tP_1(\partial^2) + 
\tP_2(\frac{x+y}{x-y} u^{-1} \Delta (\partial_x - \partial_y)).$$

Recall this means
\begin{align} \label{key_comp}
   [P(y^2), P(hx^m) ] &=m^2P(hx^m) +  2mP(hx^my) + 
m \sum_{i<j}h_i u_i^{-1} \Delta_{ij}
\frac{x_i + x_j}{x_i - x_j} (x_i^m - x_j^m).
\end{align}

\medskip

Let $\be = \frac{1}{n!}\sum_{w \in S_n} w$ be the projector onto the
trivial representation of $S_n$, and $\be\ushecke_n\be$ be the spherical Cherednik 
algebra. We regard $ \Sym^n H_\Gamma$ and $ \Sym^n H[y]$ as subalgebras
of  $\be\ushecke_n\be$
via the maps $f \mapsto f \be$.
$\be\ushecke_n\be$ inherits the filtration by order of differential
operator from $\ushecke$. (This is 
the filtration induced by regarding it as a subalgebra of the differential
part of $\End \Sym^n H_\Gamma$.) Then corollary  \ref{diff_lem} implies
$$ \Gr \be\ushecke_n\be = \Sym^n H[x,x^{-1},p ],$$
and the proof of lemma \ref{dif_lem_0} shows that
$ \Gr \be\ushecke_n\be$ is isomorphic to 
$(\Gr \Diff \End H_\Gamma^{\otimes n})^{S_n}$
as Poisson algebras.
In other words
the symbol of $P(hx^my^b)$ is the 
symbol of $P(hx^m\partial^b)$,
and hence
\begin{align}\label{poiss}
\{P(hx^ap^b), P(x^cp^d)\} = 
-(ad-bc)P(hx^{a+c}p^{b+d-1}).
\end{align}

\medskip

Let $\CS \subseteq \Diff\End\Sym^n H_\Gamma $
be the subspace spanned by $P(y^2)$ and $P(hy)$, for $h \in H$, and 
set $\oCS = \langle \Sym^n H_\Gamma,\CS \rangle
\subseteq \Diff\End\Sym^n H_\Gamma $.
As the third term of \ref{key_comp} is clearly in $\Sym^n H_\Gamma$,
we have $P(hx^my) \in \oCS$ for all $m \in \Z$, $h \in H$. In fact:

\begin{proposition} \label{spherical bit} 
 $\oCS 
= \langle \Sym^n H_\Gamma, \Sym^n H[y] \rangle = 
\be\ushecke_n\be $.
\end{proposition}
\begin{proof}
It suffices to show that $\Gr \oCS$ is $\Gr \be\ushecke_n\be$.
As $\Gr \oCS$ is a Poisson subalgebra of  $\Sym^n H[x,x^{-1},p ]$
which contains $P(hx^mp)$ and $P(hx^m)$ for all $m \in \Z$, $h \in H$,
it suffices to show that these generate $\Sym^n H[x,x^{-1},p ]$ as a 
Poisson algebra. But \ref{poiss} implies $P(hx^ap^b)$ is in $\Gr \oCS$
for all $a$, $b$, and lemma \ref{surj} then shows that the commutative
algebra generated by this is all of  $\Sym^n H[x,x^{-1},p ]$.
\end{proof}
{\it Variant:} All the above is still true if we take $\CS$ to
be the subspace spanned by $P(y)$ and $P(hy^2)$, for $h \in H$.

\medskip

Note that the inclusions
$(\Diff \End H_\Gamma^{\otimes n})^{S_n} \subset \Diff\End \Sym^n H_\Gamma$,
$\be\ushecke_n\be
  \subset \Diff_H\End \Sym^n H_\Gamma$
are proper (for example, the latter contains the Opdam shift operators).

\begin{proposition} The algebra of operators local with respect to $\CS$ 
is $\Sym^n H[y]$.
Moreover $$ Z_{\Diff_H\End\Sym^n H_\Gamma}(\CS) =  \Sym^n H[y].$$
\end{proposition}
\begin{proof}
Let $R = \Diff \End_{H^{\otimes n}}H_\Gamma^{\otimes n}[\frac{1}{\disc}]$,
so $\Gr R = (H[x^{\pm 1},p])^{\otimes n}[\frac{1}{\disc}]$. 
We first show that the Poisson centralizer in $\Gr R$ 
of $\sum p_i^2$ is $H_\Gamma^{\otimes n}[p_1,\dots,p_n]$. 
Suppose $s \in \Gr R$, and $\{ \sum p_i^2, s \} = 0$. 
Write $s = a(x,p)/b(x)$, where 
$b(x) = \prod_{i < j} (1 -x_ix_j^{-1})^{m_{ij}}$,
$a(x) \in H_\Gamma[p]^{\otimes n}$,
and $a(x,p)$ and $b(x)$ have no common divisor.
Then
$$ 0 = \Big\{ \sum p_i^2, s  \Big\} = 
b^{-2}\Big(b \sum p_i \partial_i a - a \sum p_i \partial_i b\Big). $$
Now, for $\gamma \in \Gr R$, $\partial_i(\gamma) = 0$ for all $i$ if and only
if $ \gamma \in  (H_\Gamma[p])^{\otimes n}$. 

Hence if $\partial_i b = 0$ for all $i$, we get $\sum p_i (\partial_i a) = 0$,
so $\partial_i a = 0$ for all $i$, and $a(x,p) \in  (H_\Gamma[p])^{\otimes n}$
as claimed.
Suppose there is an $i$ with $\partial_i b \neq 0$. Then $b$ divides
$ a \sum p_i \partial_i b$, whence (as $a$ and $b$ have no common factor)
$b$ divides $\partial_ib$. This holds even though $H_\Gamma^{\otimes n}$
is not a UFD, as $b$, $\sum p_i \partial_i b \in \C[x^{\pm 1},p]^{\otimes n}$.

But for any $l \in \Z^n$,
$$ \partial_l b = \big(\sum_{\alpha > 0} m_\alpha \langle - \alpha, l 
\rangle\big) b 
+ \sum_{\alpha > 0} m_\alpha \langle  \alpha, l \rangle
\frac{b}{1 - e^{-\alpha}}$$
so that if $m_{ij} \neq 0$, $(1-x_ix_j^{-1})^{m_{ij}}$ does
not divide $\partial_i b$; a contradiction.

\smallskip

Now let $D \in \Diff^s_H\End\Sym^n H_\Gamma$ and suppose $[\CS,D] = 0$.
Then $\{ \sum p_i^2, \sigma D \} = 0$, where $\sigma D$ is the image of $D$
in $\Gr R$. Hence $\sigma D = a(p) \in H_\Gamma[p]^{\otimes n}$.
As $D$ preserves $\Sym^n H_\Gamma$, $\{ \sigma D, \cdot \}$ does also,
and so $a(p) \in \Sym^n H_\Gamma[p]$. It follows that 
$D - a(y) \in \Diff^{s-1}_H\End\Sym^n H_\Gamma$ and  $[\CS,D- a(y)] = 0$.
Hence the centralizer of $\CS$ is $\Sym^n H[y]$. Finally,
just observe that the previous proposition implies that $\Sym^n H[y]$ is local
with respect to $\CS$.
\end{proof}

\subsection{Stabilisation}
\label{sec_stable}
Let $\CS$ be the subspace of $\Sym H[y]$ spanned by $P(y^2)$ and $P(hy)$, for $h \in H$.
The results of the previous sections define a map
$$ P : \Sym H[y] = \Sym^\infty H[y]_+ \to \pdiff \Sym H_\Gamma, $$
and show
$$ Z_{\pdiff \Sym H_\Gamma}(\CS) = \Sym H[y] = \Loc_{\pdiff \Sym H_\Gamma}(\CS). $$
(To see this, observe that to get stabilisation maps we must work in the category $\alg_N$
of non-unital algebras. Equivalently, we use the functor ${}_+$ to translate this back into the world
of augmented unital algebras.)

Now  take the Dunkl-Cherednik operators with parameter $u^2$, 
i.e.{} $y_i = y_i(u^2)$,
to get a description of the Calogero-Sutherland operator in terms of the Dunkl-Cherednik
operators:
$$
\Phi_u^{-1} \, \lehn(H, u - eu^{-1}) \, \Phi_u =
 P(uy^2
) - \langle \rhotw(u^2), \rhotw(u^2) \rangle.
$$
Note that this makes sense: 
$  \langle \rhotw(u), \rhotw(u) \rangle = 2u^{-2} ( \sum_{i<j} \Delta_{ij}^2 
+ \sum_{i < j < k} \Delta_{ij} \Delta_{jk}) \in F^3\Sym H_\Gamma$.

We summarize the conclusions of the previous section:
\begin{corollary} 
Fix $u \in H^\times$, and put $K = u - eu^{-1}$. Then
$$\IM_{K} = Z_{\pdiff_{H,K} \Sym H_\Gamma}\big(\lehn, P(\partial)\big) = 
\Loc_{\pdiff_{H,K} \Sym H_\Gamma}\big(\lehn, P(\partial)\big) =
 \Phi_u \Sym H[y(u^2)] \Phi_u^{-1}.$$
Moreover  
$\IM_K \subseteq \Loc_{\Diff \End\fock} \big(\lehn, P(\partial)\big)$,
and
$\Gr \IM_K = \Sym H[p].$
\end{corollary}
Notice that $\IM_K$, defined as operators on 
$\Sym H_\Gamma$,  preserves $\fock$.
\begin{proof} 
A variant 
of proposition \ref{spherical bit} shows $\IM_K$ is local in $\Diff\End \fock$.
To see $\Gr \IM_K = \Sym H[p]$, observe  $\Gr H[y(u^2)] = \Sym H[p]$, and 
$\Phi_u \Sym H[p]\, \Phi_u^{-1} = \Sym H[p]$. The rest has been proved.

\end{proof}

Now let $(K_\lambda,u_\lambda)$ be a degeneration direction.
\begin{corollary} \label{whatever}
Let  
$\IM_K = \lim_{\lambda \to 0}\Phi_{u_\lambda} \Sym H[y(u_\lambda^2)] 
\Phi_{u_\lambda}^{-1}$.
Then 
$\IM_K \subseteq 
\Loc_{\pdiff_{H,K} \Sym H_\Gamma} \big(\lehn, P(\partial)\big), $ and
$\IM_K \subseteq 
\Loc_{\Diff \End\fock} \big(\lehn, P(\partial)\big)$.
\end{corollary}

\section{CS-locality on the Fock space $\fock$}

The definition of {\it locality} with respect to $\lehn$ (section
\ref{loc_sec}) encodes the notion of differential operators which can
be built out of functions and $\lehn$.  In this section we show that
the symbols of local operators on $\fock$ are contained in $\Sym
H[p]$; then we show that the Hecke algebra construction produces
enough (all) the  local operators.

We retain the notation of the previous sections, so 
$H$ denotes a weak Frobenius algebra, $K \in H$,
$\fock = \fock(H) = \Sym H[x]/H \Sym H[x]$, and we set
$\CS_K = \langle \lehn, P(h\partial) \mid h \in H\rangle$,
$\oCS_K = \langle \fock, \CS_K \rangle =  \langle P(x), \CS_K \rangle$
(see proposition \ref{spherical bit}).

\medskip

Filter $\Diff\End \fock$ by order of differential operator, so
that $\Gr \Diff\End \fock$ is a Poisson algebra with bracket denoted
$\{,\}$, and $\Gr^0\Diff\End \fock = \fock$. Note that 
 $\Gr^i\Diff\End \fock$ is of uncountable dimension for $ i > 0$.
Write 
$P(hx^np^m)$ for the symbol of $P(hx^n\partial^m)$.

Observe that $\Gr \oCS_K  = \langle P(x),
P(hp), P(p^2) \mid h \in H \rangle \supseteq \Sym xH[x]$. 
This implies part (ii)
of the following lemma; part (i) is immediate from the definitions.

\begin{lemma}  \label{obvious}
i) If $s \in  \Gr^n\Diff\End \fock$,
and $\{s, P(hx^m) \} = 0$ for all $h \in H$, $m \geq 0$, then
$s \in \fock$.  

If in addition $\{s,P(p)\} = 0$ then $s = 0$. If $n \neq 0$, then $s =0$.

ii) If $s, s' \in \Gr^n\Diff\End \fock$, and 
$\{s-s', P(x) \} = \{s-s', P(hp) \} = \{s-s',P(p^2)\} = 0$,
then $s = s'$.
\end{lemma}

\begin{proposition} \label{estimate}
 $ \Gr  \Loc_{\Diff End  \fock}(\CS_K) \subseteq \Sym H[p]$.
\end{proposition}
\begin{proof} Inductively define $L^n \subseteq \Gr^n\Diff\End \fock$
by setting $L^{-1} = 0$, and 
$$L^n = \{ s \in \Gr^n\Diff\End \fock \mid
\{s, P(p^2) \} = \{s, P(p) \} = 0, \text{ and }
\{s, P(x) \} \in \oplus_i \Gr^i \oCS_K . L^{n-i-1} \}.$$
Obviously $ \Gr  \Loc_{\Diff End  \fock}(\CS_K) \subseteq L$,
so it suffices to
show that if $L^{r} \subseteq \Sym H[p]$, then 
 $L^{r+1} \subseteq \Sym H[p]$.

Write $\Gr^{r}\oCS[a]$ for the eigenspace of $\{P(p),\cdot\}$
with eigenvalue $a$; for example  $\Gr^{r}\oCS[1]$ is spanned by monomials
$P(hxp^N)P(h_1p^{n_1})\dots P(h_kp^{n_k})$ with $N + \sum n_i = r$.

If $s\in L^{r+1}$, then $\{s, P(x) \} \in \Gr^{r}\oCS$
as $L^{r} \subseteq \Gr^{r}\oCS$ by inductive assumption.
Moreover, as $\{s,P(p) \} = 0$,  $\{s, P(x) \} \in \Gr^{r}\oCS[1]$.
As
$$\{\cdot, P(x) \} : \Gr^{r+1} \Sym H[\partial] \to \Gr^{r}\oCS[1] $$
is surjective, there exists a $\phi \in \Gr^{r+1} \Sym H[\partial]
 = (\Sym H[p])^{r+1}$
with $\{s- \phi, P(x) \} = 0$. Clearly 
$\{s - \phi, P(p) \} =  \{s - \phi, P(p^2) \} = 0$.
The previous lemma implies $s = \phi$.
\end{proof}

\begin{theorem} \label{construct}
$\IM_K =  \Loc_{\Diff End  \fock}(\CS_K)$. In particular, 
$\IM_K$ is independent of the choice of degeneration direction.
Moreover $\Gr \IM_K = \Sym H[p]$.
\end{theorem}
\begin{proof}
Corollary \ref{whatever} shows that
$$\IM_K 
\subseteq \Loc_{\Diff End  \fock}(\CS_K).$$
But the previous proposition states
$$ \Gr  \Loc_{\Diff End  \fock}(\CS_K) \subseteq \Sym H[p], $$
so $\Gr \IM_K \subseteq  H[p]$.
As the degree of a differential operator can only decrease when we degenerate,
$\dim {\IM_K}^{\leq s} \geq 
\dim {\IM_{u-eu^{-1}}}^{ \leq s} = \dim H[p]^{\leq s}$, 
where $X^{\leq s}$ denotes differential operators of degree $\leq s$.

Hence we have $\Gr \IM_K = H[p]$, and the inclusions above are equalities.

\end{proof}

\section{Hilbert schemes of points on surfaces.}
\label{section_hs}
Let $X$ be a smooth, connected projective surface over $\C$. 
We recall the basic results about the geometry and
cohomology of the Hilbert scheme of points on $X$.

\begin{definition} \label{def_hilbert_scheme}
Denote by $\Hilb n X$ the moduli space of zero-dimensional
subschemes $Z \subset X$ satisfying
\begin{equation*}
\dim H^0(\Oo_Z) = n
\end{equation*}
\end{definition}
$\Hilb n X$ is a projective variety if $X$ is.
What is special about the geometry of Hilbert schemes of points on
surfaces are the following theorems:
\begin{theorem}[Fogarty \cite{fogarty}]
$\Hilb n X$ is a smooth, connected, irreducible, projective
variety.
\end{theorem}
Let $S^n X$ the symmetric product $X^n / S_n$. There is a natural
Hilbert-Chow map
\begin{equation*}
\Hilb n X \to S^n X
\end{equation*}
which sends a subscheme $Z \subset X$ to it's support, counted
with multiplicity. This map is a resolution of singularities;
its main property is:
\begin{theorem}[Briancon, \cite{bria}]
If $x \in X$, the fiber of the Hilbert-Chow morphism over $nx$
is an irreducible variety of dimension $n-1$.
\end{theorem}
Both these theorems are false if $\dim X > 2$.

Since $\Hilb n X$ represents a functor, there is a universal
family of subschemes
\begin{equation*}
\Xi_n \subset X \times \Hilb n X
\end{equation*}
flat and finite over $\Hilb n X$, of degree $n$. For each vector
bundle $V$ on $X$, $\Xi_n$ induces the tautological vector bundle
\begin{equation*}
V^{[n]}  = p_{2\ast}p_1^{\ast} \Oo_{\Xi_n}
\end{equation*}
on $\Hilb n X$, with $\rank V^{[n]} = n \rank V$.

\subsection{Fock space structure on the Hilbert scheme}
We recall certain algebraic structures on the
cohomology of the Hilbert schemes induced from correspondences.
Let
\begin{equation*}
\H = \bigoplus_{n, i} H^{i} (\Hilb n X, \C).
\end{equation*}
This is a super vector space, graded by cohomology degree.
It is also graded by length of subscheme: 
set $\H^n = \oplus_i H^{i} (\Hilb n X, \C)$.
The bilinear form on $H^{\ast}(\Hilb n X)$ induces a
non-degenerate, super-symmetric bilinear form on $\H$.

Define
\begin{equation*}
E_{n,m}^{n+m} \subset \Hilb n X \times \Hilb m X \times \Hilb
{n+m} X
\end{equation*}
to be the closure of the locus of triples $(a,b,c)$ of subschemes
of $X$, satisfying:
\begin{align*}
a \cap b &= \emptyset\\
a \cup b &= c
\end{align*}
Since $E_{n,m}^{n+m}$ is a closed subscheme, it has a fundamental
class
\begin{equation*}
[E_{n,m}^{n+m}] \in H^{\ast}(\Hilb n X , \C) \otimes
H^{\ast}(\Hilb m X , \C) \otimes H^{\ast}(\Hilb {n+m} X , \C)
\end{equation*}
Pulling back, intersecting with $[E_{n,m}^{n+m}]$ and pushing
forward induces maps
\begin{align*}
m : H^{\ast}(\Hilb n X , \C) \otimes H^{\ast}(\Hilb m X , \C)
    &\to H^{\ast}(\Hilb {n+m} X , \C) \\
c : H^{\ast}(\Hilb {n+m} X , \C)
    &\to H^{\ast}(\Hilb n X , \C) \otimes H^{\ast}(\Hilb m X , \C)
\end{align*}
Putting all these together, we get maps
\begin{align*}
m : \H \otimes \H &\to \H \\
c : \H &\to \H \otimes \H
\end{align*}

These maps are adjoint : $m^{\dagger} = c$.
\begin{theorem}[Grojnowski \cite{groj}]
$m$ and $c$ give $\H$ the structure of a commutative,
cocommutative graded Hopf algebra, in the super sense. The unit is
given by the identity in the ring $H^{\ast} (\Hilb 0 X, \C) \iso
\C$.
\end{theorem}
An equivalent formulation 
was obtained in \cite{nak},
which describes the structure of $\H$ as a module over the Heisenberg algebra.

A commutative, cocommutative graded Hopf algebra
is naturally isomorphic to the free commutative algebra on the
space of primitive elements, 
\begin{equation*}
\Prim_n = \left\{ \alpha \in \H^n \mid c(\alpha)
    = \alpha \otimes 1 + 1 \otimes \alpha \right\}
\end{equation*}
One can identify
$\Prim_{n,\ast}$ explicitly as the image of a map induced from a
certain correspondence.  Let
\begin{equation*}
Z_n \subset X \times \Hilb n X
\end{equation*}
be the locus of pairs $(x, a)$ satisfying
\begin{equation*}
\Supp a = nx
\end{equation*}
$Z_n$ is a closed subscheme, and so induces a linear map
\begin{equation*}
[Z_n] : H^{\ast} (X , \C) \to H^{\ast} (\Hilb n X, \C) = \H^n
\end{equation*}
One can show that this map is of cohomology degree $2n-2$.
Geometrically, this map takes a cycle $A \subset X$ to the cycle
in $\Hilb n X$, consisting of those subschemes $a$ which are
supported on exactly one point, which lies in $A$.
\begin{theorem}[Grojnowski \cite{groj}, Nakajima \cite{nak}]
$[Z_n]$ induces an isomorphism
$$ H^{\ast}(X, \C) \iso \Prim_{n}(\H) $$
and hence a canonical isomorphism of graded Hopf algebras
$$ \H \iso \fock(H^{\ast}(X , \C)) = \Sym^{\ast}(  H^{\ast}(X, \C)\otimes x\C[x] ) $$

\end{theorem}
In this
isomorphism  $x$ is given degree 1, so 
this identifies $\H^n$ with the energy weight space of weight $n$.

If we take $X$ smooth, but not necessarily projective, we still have a canonical
identification $\H \iso \fock(H^{\ast}(X , \C))$. However, now this isomorphism
is only as graded coalgebras. Dually, if we intepret $H^*( , \C)$ as compactly supported 
cohomology, we again have such an isomorphism; but now it is of graded algebras.


\subsection{Lehn's theorem}
Cup product on $H^{\ast}(\Hilb n X, \C) = \fock^n$ induces a ring structure
$$
\star: \fock^n \otimes \fock^n \to \fock^n.
$$
This is distinct from the Fock space structure induced by correspondances. 
We would like to describe it.

Let $T^*_X $ be the cotangent bundle of $X$, and $T^{\ast [n]}_X$ be the 
associated tautological bundle on $\Hilb n X$, a bundle of rank $2n$.
Write $H = H^{\ast}(X, \C)$, and $K = c_1(T_X^*)$ for the canonical class.
$H$ is a weak Frobenius algebra; if $X$ is projective it is even a Frobenius algebra.
Note moreover that $e = c_2(T_X^*)$, and so an element $u$ such that
$u^2 - Ku - e = 0$ is precisely a Chern root of $T^*_X$. We define the 
Calogero-Sutherland operator $\lehn(H,K) \in \End \fock$. The main result of \cite{lehn}
is the computation of cup product with the boundary of the Hilbert scheme in terms of 
the Fock space coordinates of \cite{groj, nak}: 

\begin{theorem}[Lehn \cite{lehn}] \label{lehn_th}
The linear map $\fock \to \fock$ defined on $\fock^n$ as 
$$ x \mapsto c_1( T^{\ast[n]}_X )  \star x $$
is the Calogero-Sutherland operator $\lehn(H, K)$.
\label{lehnop}
\end{theorem}
This is theorems 3.10 and 4.2 of \cite{lehn}.

The following theorem is due to Lehn \cite{lehn} for the subalgebra
$H^\ast_{alg} (X,\C) \subset H^\ast(X,\C)$ of algebraic cohomology, and was extended
to all of $H^\ast(X,\C)$ by Li, Qin and Wang \cite{lqw, lqw2}. 
Define $ \CS_K = \langle \lehn, P(h\partial) \mid h \in H\rangle$.

\begin{theorem} \label{chern class operators}
There are linear maps 
\begin{equation*}
\cht_i : H^{\ast}(X, \C) \to \Diff^{i+1} \End \fock
\end{equation*}
for $i \ge 0$, such that  $\cht_i(h)$ preserves each subspace $\fock_n$, and
which satisfy 
\begin{equation*}
[\cht_i(h), \lehn ]  = 0, \qquad \cht_i(h) (1) = 0, \quad  \text{ and } \quad 
[\cht_i(h), P(h' x) ] = (\Ad \lehn)^i (P(hh' x)).
\end{equation*}
Write $\dhilb$ for the subalgebra of $\Diff \End \fock$ generated by 
$\cht_i(h)$ for $h \in H$, $i \geq 0$. 
Then
 $\CS_K \subseteq \dhilb$,
and 
the image of $\dhilb$ in $\End \fock^n$ coincides with 
the algebra of left multiplication
operators, i.e.{} the 
image of the natural map
$H(X,\C) \to \End \fock^n$, $ a \mapsto a*$.
\label{chern}
\end{theorem}


\section{Hecke algebras and Hilbert schemes}
We now apply the results of sections 2--5 in the context of section \ref{section_hs}.
This gives a description of the ring structure on $H^*(\Hilb n X, \C)$ in the 
Fock space coordinates.
In particular, this gives an explicit algebraic construction
of  the ring $\dhilb$, and shows it 
depends only on the weak Frobenius algebra $H = H(X,\C)$ and the class $K \in H$.

\begin{theorem}  \label{loc_thm}
 Let $X$ be a smooth algebraic surface, $H = H(X,\C)$ its cohomology ring.
Identify $\fock^n \iso H^*(\Hilb n X, \C)$ as above. 
Let $\IM_K$ be the integrals of motion;
recall that this is independent of the choice of degeneration direction.

Then 

$$ i) \qquad
\dhilb = \Loc_{\Diff End  \fock}(\CS_K), \text{ and } $$

$$ ii) \qquad
\Loc_{\Diff End  \fock}(\CS_K) = \IM_K. $$ 

In particular, 
the image of $\IM_K$ in $\End \fock^n$ coincides with 
the algebra of left multiplication
operators, i.e.{} the 
image of the natural map
$H(X,\C) \to \End \fock^n$, $ a \mapsto a*$.
\end{theorem}
\begin{proof}
We have already proved (ii) holds for any weak Frobenius algebra $H$,
and $K \in H$ as theorem \ref{construct}. We prove (i).

By theorem \ref{chern}, $\CS_K \subseteq \dhilb \subseteq 
Z_{\Diff \End \fock}(\CS_K)$, and $[\cht_i(h), P(x) ] \in 
\langle \CS_K , \fock \rangle$.
Moreover, as
$[[\lehn,P(x)], P(hx^n)] = 2n P(hx^{n+1})$,
$\langle \CS_K , \fock \rangle = \langle \CS_K , P(x) \rangle$.
Hence $\dhilb . \langle \CS_K , \fock \rangle$
is a subalgebra, and so $\dhilb \subseteq \Loc(\CS_K)$.

Let us compute the symbol $s$ of $\cht_i(h)$. The conditions in 
theorem \ref{chern} give
$$\{s, P(h'x) \} = \Ad P(p^2)^i P(hh'x) = 2^iP(hh'xp^i),$$
and $\{s,P(h'p) \} = \{s, P(h'p^2) \} = 0$.
By lemma \ref{obvious}(i),
$s = \frac{2^i}{i+1}P(hp^{i+1})$.

So $\Gr \dhilb = \Sym H[p]$. But $\Gr \Loc(\CS_K)  \subseteq \Sym H[p]$,
by \ref{estimate}. We must have equality.
\end{proof}

\begin{remark} The algebra $\Sym H[p]$ is Poisson self-centralising in 
$\Gr\Diff\End\fock$; it follows that $\dhilb$ is also the algebra of 
all differential operators on $\fock$ which preserve each $\fock^n$,
and whose restriction to each $\fock^n$ is cup product with {\it some}
class. (As this algebra is obviously commutative, and  contains
$\langle \cht_i \rangle$.)
\end{remark}

In this paper we have set up a formalism precisely linking the combinatorics 
of intersection theory on the Hilbert scheme with integrable systems
and certain generalizations of Cherednik algebras.

It is now straightforward to compute all geometric
information on the Hilbert scheme in terms of the well understood
combinatorics of Jack polynomials and/or  Cherednik algebras.  
Elaborations and further
generalizations will appear elsewhere.


\end{document}